\documentclass{article}
\usepackage{graphicx}
\usepackage{amsmath}

\input{tcilatex}

\begin{document}

\ \ \ \ \ \ \ \ \textbf{SELF-SIMILAR FOCUSING IN POROUS MEDIA: }

\ \ \ \ \ \ \ \ \ \ \ \ \ \ \ \ \ \ \ \ \ \ \textbf{AN EXPLICIT} \textbf{%
CALCULATION}

\ \ \ \ \ \ \ \ \ \ \ \ \ \ \ \ \ \ \ \ \ \ \ \ \ \ \ \ \ \ \ \ \ \ \ \ \ \
\ D. G. Aronson

Institute for Mathematics and its Applications, Minneapolis MN 55455

\bigskip

ABSTRACT. We consider a porous medium flow in which the gas is initially
distributed in the exterior of an empty region (a hole) and study the final
stage of the hole-filling process. From general theory it is known that this
hole-filling is asymptotically described by a self-similar solution which
depends on a constant determined by the initial configuration. In general,
this constant must be found either experimentally or numerically. Here we
give an example of a one-dimensional symmetric flow where the appropriate
constant is obtained explicitly.

PACS number(s): 47.56+r

\bigskip

Self-similarity arises in a natural way in the description of critical
behavior in various physical problems. Many examples can be found in the
book of G. I. Barenblatt [1]. In this note we consider the focusing problem
for porous medium flow. In this problem the material is initially
distributed in the exterior of an empty region (a hole) and one is
interested in the details of the final stage of the hole-filling process.
For axially symmetric flows it is known [2] that this process is ultimately
self-similar. Here we present an explicit calculation of the development of
self-similarity in the focusing of a particular one-dimensional flow.

The evolution of the (scaled) density $U(r,t)$ of an ideal gas flowing in an
axially symmetric homogeneous porous medium is governed by the degenerate
nonlinear diffusion equation 
\begin{equation*}
\partial _{t}U=(\partial _{rr}^{2}+\frac{d-1}{r}\partial _{r})(U^{m}),
\end{equation*}
where $d$ is the dimension of the space and $m$ is a constant derived from
the ideal gas law. Using the ideal gas law we can replace the density $U$ by
the (scaled) pressure 
\begin{equation*}
V=\frac{m}{m-1}U^{m-1}
\end{equation*}
whose evolution is governed by the equation 
\begin{equation}
\partial _{t}V=(m-1)V(\partial _{rr}^{2}+\frac{d-1}{r}\partial
_{r})V+(\partial _{r}V)^{2}.  \tag{1}
\end{equation}
In the context of gas flow $m>2$, but other values of $m$ occur in various
applications. Here we will assume that $m>1$ which assures that the speed of
propagation of disturbances from rest is finite (slow diffusion). For $m\leq
1$ the corresponding speed is infinite (fast diffusion) [3].

In the focusing problem for (1) we seek a solution in $[0,\infty )\times
(\tau ,\infty )$ for some $\tau \in \mathbf{R}$ such that 
\begin{equation}
V(\cdot ,\tau )=V_{0}(\cdot )\text{ on }[0,\infty ),  \tag{2}
\end{equation}
where $V_{0}$ is a given function which satisfies 
\begin{equation*}
V_{0}(r)\left\{ 
\begin{tabular}{l}
$=0$ on $[0,a]\cup \lbrack b,\infty )$ \\ 
$>0$ on $(a,b)$%
\end{tabular}
\right.
\end{equation*}
for some $0<a<b<\infty .$ It is known [3] that problem (1),(2) possesses a
unique continuous generalized solution. As time increases from $t=\tau $
material flows outward from $r=b$ \ and inward from $r=a$. There is a
non-increasing inner interface curve $r=a(t)$ and a non-decreasing outer
interface curve $r=b(t)$ with $a(\tau )=a,b(\tau )=b,$ and 
\begin{equation*}
V(r,t)\left\{ 
\begin{tabular}{l}
$=0$ on $[0,a(t)]\cup \lbrack b(t),\infty )$ \\ 
$>0$ on $(a(t),b(t))$%
\end{tabular}
\right. .
\end{equation*}
Both interface curves become monotone in finite time and there is a finite $%
T>\tau $ such that $a(t)>0$ for $t<T$ and $a(T)=0.$

For normalization we will assume that $T=0$ so that $\tau <0.$ It is shown
in [4] that there exists a one-parameter family $\{g_{c}(r,t)\}$ of
self-similar solutions to (1) defined for $c\in \mathbf{R}^{+}$ and $%
(r,t)\in \lbrack 0,\infty )\times (-\infty ,0]$ which focus at $t=0$ (the
Graveleau solutions). Specifically, there exist numbers $\alpha ^{\ast
}(d,m)\in \lbrack 1,2)$ and $\gamma (d,m)\in \mathbf{R}^{-}$ such that for $%
t<0$%
\begin{equation}
g_{c}(r,t)=\frac{r^{2}}{-t}\varphi (c\eta ),  \tag{3}
\end{equation}
where 
\begin{equation*}
\varphi \left\{ 
\begin{tabular}{l}
$>0$ on $(\gamma ,0)$ \\ 
$=0$ on $(-\infty ,\gamma )$%
\end{tabular}
\right.
\end{equation*}
and 
\begin{equation*}
\eta =tr^{-\alpha ^{\ast }}.
\end{equation*}
$g_{c}$ satisfies 
\begin{equation*}
g_{c}(r,t)\left\{ 
\begin{tabular}{l}
$=0$ for $r\in \lbrack 0,\rho _{c}(t)]$ \\ 
$>0$ for $r\in (\rho _{c}(t),\infty )$%
\end{tabular}
\right. ,
\end{equation*}
where 
\begin{equation*}
x=\rho _{c}(t)\equiv \left( \frac{ct}{\gamma }\right) ^{1/\alpha ^{\ast }}
\end{equation*}
is the interface. The function $\varphi $ and the similarity exponent $%
\alpha ^{\ast }$ are obtained by solving a nonlinear eigenvalue problem and,
in general, must be found numerically. However, in one space dimension $%
(d=1) $ we have $\alpha ^{\ast }(1,m)=-\gamma (1,m)=1$ and each Graveleau
solution is a pair of converging plane waves given by 
\begin{equation}
g_{c}(\left| x\right| ,t)=c\{\left| x\right| +ct\}_{+},  \tag{4}
\end{equation}
where $\{\cdot \}_{+}=\max (0,\cdot )$. Thus, for $d=1$, 
\begin{equation}
\eta =\frac{t}{\left| x\right| }\text{ and }\varphi (\zeta )=-\zeta (1+\zeta
).  \tag{5}
\end{equation}

The Graveleau solutions describe the asymptotics of the focusing problem
(1),(2) in the following sense [2].

\textit{Let }$V$\textit{\ denote the generalized solution of the initial
value problem (1),(2). There exists a }$c^{\ast }\in \mathbf{R}^{+}$ (%
\textit{depending only on} $d,m,$\textit{\ and }$V_{0}$) \textit{such that
the inner interface } 
\begin{equation*}
a(t)\sim \rho _{c^{\ast }}(t)\text{\textit{as }}t\nearrow 0.
\end{equation*}
\textit{For each fixed }$\eta \in (-\infty ,0]$\textit{\ the profile} $V$%
\textit{\ approaches the profile }$g_{c^{\ast }}$\textit{\ as }$r\searrow 0$%
, \textit{i.e.,} 
\begin{equation}
V(r,\eta r^{\alpha ^{\ast }})\sim g_{c^{\ast }}(r,\eta r^{\alpha ^{\ast }})%
\text{\textit{as }}r\searrow 0.  \tag{6}
\end{equation}

The constant $c^{\ast }$ depends on the initial function\ $V_{0}$ and so
must be determined on a case-by-case basis, usually numerically. However, in
one particular case we can determine $c^{\ast }$explicitly. Consider the
focusing problem for $d=1$ where we start at some time $\tau <0$ with
symmetrically placed point masses at $x=\pm \xi $ and assume that focusing
takes place at $T=0.$ In view of the symmetry it will suffice to consider $%
x>0.$

The solution to (1) corresponding to a mass $M$ initially concentrated at $%
(\xi ,\tau )$ is given by 
\begin{equation*}
V_{M}(x,\xi ,t,\tau )=\frac{\beta R(t-\tau )^{2}}{2(t-\tau )}\left\{ 1-\frac{%
(x-\xi )^{2}}{R(t-\tau )^{2}}\right\} _{+},
\end{equation*}
where 
\begin{eqnarray*}
\beta &=&\frac{1}{m-1},B=\frac{m-1}{2m(m+1)},R(t)=\sqrt{\frac{A}{B}}t^{\beta
}, \\
M &=&A^{(m+1)/2(m-1)}B^{-1/2}\int_{0}^{\pi /2}(\cos \theta
)^{(m+1)/(m-1)}d\theta .
\end{eqnarray*}
Note that $V_{M}(x,\xi ,t,\tau )>0$ only on the interval 
\begin{equation}
\xi -R(t-\tau )<x<\xi +R(t-\tau ).  \tag{7}
\end{equation}
This solution is due to Barenblatt [5] and is a self-similar solution of the
first kind meaning essentially that the similarity exponent can be found by
dimensional analysis. In contrast the Graveleau solution is a self-similar
solution of the second kind whose similarity exponent $\alpha ^{\ast }$
cannot be obtained by dimensional analysis alone [1].

In the focusing problem we require that the left hand boundary of the
support of $V_{M}$ just reaches the origin at $t=0,$ i.e., that 
\begin{equation}
\xi =R(-\tau )=\sqrt{\frac{A}{B}}(-\tau )^{\beta }.  \tag{8}
\end{equation}
For given $(\xi ,\tau )$ this condition determines the admissable mass since 
\begin{equation*}
A=\frac{B\xi ^{2}}{(-\tau )^{2\beta }}.
\end{equation*}
Using (8) we can rewrite $V_{M}$ as 
\begin{equation}
V_{M}(x,\xi ,t,\tau )=\frac{\xi ^{2}\beta }{2(-\tau )}\left( 1-\frac{t}{\tau 
}\right) ^{2\beta -1}\left\{ 1-\frac{(1-x/\xi )^{2}}{(1-t/\tau )^{2\beta }}%
\right\} _{+}.  \tag{9}
\end{equation}
Note that, in view of (7), $V_{M}(x-\xi ,t-\tau )\rightarrow 0$ as $%
x\rightarrow 0$ for any $t\in (\tau ,0].$

If we introduce the dimensionless quantities 
\begin{equation*}
\Pi =\frac{(-\tau )V_{M}}{\beta \xi ^{2}},\Pi _{1}=\frac{t}{\tau },\Pi _{2}=%
\frac{x}{\xi }
\end{equation*}
we can rewrite (9) in the form 
\begin{equation*}
\Pi =\Phi (\Pi _{1},\Pi _{2})\equiv \frac{1}{2}(1-\Pi _{1})^{2\beta
-1}\left\{ 1-\frac{(1-\Pi _{2})^{2}}{(1-\Pi _{1})^{2\beta }}\right\} _{+},
\end{equation*}
where $\Pi \rightarrow 0$ as $\Pi _{2}\rightarrow 0$ for each $\Pi _{1}.$
Set 
\begin{equation*}
\Pi ^{\ast }=\Pi _{1}/\Pi _{2}.
\end{equation*}
Then 
\begin{equation*}
\Pi =\frac{1}{2}(1-\Pi ^{\ast }\Pi _{2})^{2\beta -1}\left\{ 1-\frac{(1-\Pi
_{2})^{2}}{(1-\Pi ^{\ast }\Pi _{2})^{2\beta }}\right\} _{+}.
\end{equation*}
With $\Pi ^{\ast }$ held fixed, it follows from Taylor's theorem that 
\begin{equation*}
\Pi =\Pi _{2}(1-\beta \Pi ^{\ast })+O(\Pi _{2}^{2})\text{ as }\Pi
_{2}\rightarrow 0.
\end{equation*}
Therefore 
\begin{equation}
\Phi (\Pi _{1},\Pi _{2})\sim \Pi _{2}(1-\beta \Pi ^{\ast })\text{ as }\Pi
_{2}\rightarrow 0\text{ with }\Pi ^{\ast }\text{ constant.}  \tag{10}
\end{equation}

To interpret this result in the original dimensional variables we set 
\begin{equation}
\eta =\frac{t}{x}\text{ and }p=\frac{\xi }{-\tau }.  \tag{11}
\end{equation}
Note that $\eta <0,p>0,$ and $\Pi ^{\ast }=-p\eta .$ It follows from (10)
and (11) that for each fixed $\eta $%
\begin{equation*}
V_{M}(x,\xi ,t,\tau )\sim p\beta x\left\{ 1+p\beta \eta \right\} _{+}\text{
as }x\rightarrow 0.
\end{equation*}
Using (5) we can rewrite this as 
\begin{equation*}
V_{M}(x,\xi ,t,\tau )\sim \frac{x^{2}}{-t}\varphi (p\beta \eta )=g_{p\beta
}(\left| x\right| ,t).
\end{equation*}
Thus (6) holds with 
\begin{equation*}
c^{\ast }=p\beta .
\end{equation*}
Similarly, for the inner interface we have $\Pi ^{\ast }\sim 1/\beta $ so
that 
\begin{equation*}
a(t)\sim p\beta (-t)=\rho _{p\beta }(t).
\end{equation*}

Two aspects of the explicit computation of $c^{\ast }$ are essential. Since
we are dealing with a symmetric one-dimensional problem we need only
consider $x>0$ in the pre-focusing regime. This allows us to exploit the
explicitly known Barenblatt point mass solution. A similar computation could
be carried out if there were other explicit focusing solutions.

\bigskip

[1] G. I. Barenblatt, \textit{Scaling, Self-Similarity, and Intermediate
Asymptotics,}

\ \ \ \ Cambridge University Press,1996.

[2] S. B. Angenent \& D. G. Aronson, Comm. P. D. E., \textbf{20}(1995),
1217-1240.

[3] D. G. Aronson, \textit{The Porous Medium Equation, }Lecture Notes in

\ \ \ \ Mathematics 1224, Springer Verlag, 1986.

[4] D. G. Aronson \& J. Graveleau, Euro. J. Appl. Math., \textbf{4}(1993),
65-81.

[5] G. I. Barenblatt, Prikl. Mat. Meh., \textbf{16}(1952), 67-78.

\ \ \ \ \ \ \ \ \ \ \ \ \ \ \ \ \ \ \ \ \ \ \ \ \ \ \ \ \ \ \ \ \ \ \ \ \ \
\ \ \ \ \ \ \ \ \ \ \ \ \ \ \ \ \ \ \ \ \ \ \ \ \ \ \ \ \ \ \ \ \ \ \ \ \ \
\ \ \ \ \ \ \ \ \ \ \ \ \ \ \ \ \ \ \ \ \ \ \ \ \ \ \ \ \ \ \ \ \ \ \ \ \ \
\ \ \ \ \ \ \ \ \ \ \ \ \ \ \ \ \ \ \ \ \ \ \ \ \ \ \ \ \ \ \ \ \ \ \ \ \ \
\ \ \ \ \ \ \ \ \ \ \ \ \ \ \ \ \ \ \ \ \ \ \ \ \ \ \ \ \ \ \ \ \ \ \ \ \ \
\ \ \ \ \ \ \ \ \ \ \ \ \ \ \ \ \ \ \ \ \ \ \ \ \ \ \ \ \ \ \ \ \ \ \ \ \ \
\ \ \ \ \ \ \ \ \ \ \ \ \ \ \ \ \ \ \ \ \ \ \ \ \ \ \ \ \ \ \ \ \ \ \ \ \ \
\ \ \ \ \ \ \ \ \ \ \ \ \ \ \ \ \ \ \ \ \ \ \ \ \ \ \ \ \ \ \ \ \ \ \ \ \ \
\ \ \ \ \ \ 

\end{document}